\documentclass[11pt,twoside,a4paper]{article}
\title{} \author{} \date{}
\usepackage{latexsym,amssymb,times}
\input amssym.def
\newtheorem{te}{Theorem}[section]

\newtheorem{fac}[te]{Fact}

\newtheorem{rem}[te]{Remark}

\def\dok{\noindent{\bf Proof. }}
\def\kdok{\hfill $\Box$ \par \vspace*{2mm} }
\def\a{\alpha}

\def\f{\varphi}
\def\p{\psi}
\def\o{\omega}
\def\k{\kappa}

\def\r{\rho}

\def\t{\tau}

\def\ve{\varepsilon}

\def\N{{\mathbb N}}
\def\X{{\mathbb X}}
\def\Y{{\mathbb Y}}
\def\Z{{\mathbb Z}}

\def\BI{{\mathbb I}}

\def\BJ{{\mathbb J}}

\def\CC{{\mathcal C}}
\def\CD{{\mathcal D}}
\def\CT{{\mathcal T}}
\def\I{{\mathcal I}}

\def\c{{\mathfrak{c}}}
\def\la{\langle}
\def\ra{\rangle}
\def\Aut{\mathop{\rm Aut}\nolimits}
\def\Mod{\mathop{\rm Mod}\nolimits}
\def\EF{\mathop{\rm EF}\nolimits}
\def\Th{\mathop{\rm Th}\nolimits}
\def\PI{\mathop{\rm PI}\nolimits}
\def\bcd{\dot{\bigcup}}
\def\du{\mathrel{\dot{\cup}}}
\def\IT{\mathop{\rm it}\nolimits}
\begin{document}
\thispagestyle{plain}
\begin{center}
          {\large \bf \uppercase{Vaught's conjecture and theories of partial order \\[2mm] admitting a finite lexicographic decomposition}}
\end{center}
\begin{center}
{\bf Milo\v s S.\ Kurili\'c}\footnote{Department of Mathematics and Informatics, Faculty of Sciences, University of Novi Sad,
                                      Trg Dositeja Obradovi\'ca 4, 21000 Novi Sad, Serbia.
                                      e-mail: milos@dmi.uns.ac.rs}
\end{center}
\begin{abstract}
\noindent
A complete theory $\CT$ of partial order is an FLD$_1$-{\it theory}
iff some (equivalently, any) of its models $\X$
admits a finite lexicographic decomposition $\X =\sum _{\BI}\X _i$,
where $\BI$ is a finite partial order and $\X _i$-s are partial orders with a largest element.
Then we write $\sum _{\BI}\X_i\in \CD (\CT)$
and call $\sum _{\BI}\X_i$ a VC-{\it decomposition} (resp.\ a VC$^\sharp$-{\it decomposition})
iff $\X _i$ satisfies Vaught's conjecture (VC)  (resp.\ VC$^\sharp$: $I(\X _i)\in \{ 1,\c\}$), for each $i\in I$.
$\CT$ is called {\it actually Vaught's}
iff for some $\sum _{\BI}\X_i\in \CD (\CT)$ there are sentences $\t _i\in \Th (\X _i)$, $i\in I$, providing VC.
We prove that:
(1) VC is true for $\CT$ iff $\CT$ is large or its atomic model has a VC decomposition;
(2) VC is true for each actually Vaught's FLD$_1$ theory;
(3) VC$^\sharp$ is true for $\CT$, if there is a VC$^\sharp$-decomposition of a model of $\CT$.
Defining FLD$_0$ theories (here $\X _i$-s have a smallest element, ``0") we obtain duals of these statements.
Consequently, since the classes
$\CC ^{\rm lo}_0\subset\CC ^{\rm tree}_0\subset \CC ^{\rm reticle}_0$  and $\CC ^{\rm ba}$
of linear orders with 0, rooted trees, reticles with 0 and Boolean algebras are first-order definable,
VC is true for the partial orders from the closure $\la \CC ^{\rm reticle}_0\cup \CC ^{\rm ba}\ra_{\Sigma}$,
where $\la \CC\ra_{\Sigma}$ denotes the closure of a class $\CC$ under finite lexicographic sums.
Defining the closure $\la \CC\ra_{\Sigma ^r}$ under finite lexicographic sums of rooted summands, $\sum _{\BI}(\X _i)_r$,
we show that $\la \CC ^{{\rm VC}^{\sharp}}\ra _{\Sigma ^r}=\CC ^{{\rm VC}^{\sharp}}$,
where $\CC ^{{\rm VC}^{\sharp}}$ is the class of all partial orders satisfying {\rm VC}$^{\sharp}$.
In particular {\rm VC}$^{\sharp}$ is true for a large class of partial orders of the form
$\sum _{\BI}(\dot{\bigcup}_{j<n_i}\prod _{k<m_i^j}\X _i^{j,k})_r$,
where $\X _i^{j,k}$-s can be linear orders,
or Boolean algebras, or belong to a wide class of trees.


\end{abstract}
\section{Introduction}\label{S1}
We recall that Vaught's conjecture (VC), stated by Robert Vaught in 1959 \cite{Vau},
is the statement that the number $I(\CT ,\o)$ of non-isomorphic countable models
of a complete countable first-order theory $\CT$
is either at most countable or continuum.
The results related to this (still open) problem include a reduction relevant to this paper:
VC is equivalent to its restriction to the theories of partial order (see \cite{Hodg}, p.\ 231).
Regarding such theories and denoting  by VC$^\sharp$ the ``sharp" version of Vaught's conjecture, $I(\CT)\in \{1,\c\}$, we recall the following classical results.
\begin{fac}\label{T721}
Vaught's conjecture is true for the theories of

(a) linear orders; moreover, VC$^\sharp$ is true (Rubin \cite{Rub});

(b) model-theoretic trees  (Steel \cite{Stee});

(c) reticles (partial orders which do not embed the four-element poset $N$) (Schmerl \cite{Sch3});

(d) Boolean algebras; moreover, VC$^\sharp$ is true (Iverson \cite{Ive}).
\end{fac}
Continuing the investigation from \cite{KMon}--\cite{Ksharp}
we consider several model-theoretic constructions (e.g.\ interpretations, direct products, etc.)
and deal with the question whether they preserve VC.
Namely, taking a class $\CC$ of structures for which VC was already confirmed,
our goal is to confirm VC for the structures from its closure $\la \CC\ra _c$ under a construction $c$.

For example, if $L$ is any relational language,
$\la \CC\ra _{\rm def}$ is the class of $L$-structures definable in structures from $\CC$ by quantifier free formulas,
$\CC^{\rm lo}$ is the class of linear orders and
$\CC^{\rm lo}_{\rm lab}$ is the class of linear orders colored into finitely many convex colors (labelled linear orders),
then by \cite{KMon}--\cite{KFMD} we have
\begin{fac}\label{T718}
VC$^\sharp$ is true for all relational structures from the class $\la \CC^{\rm lo}_{\rm lab} \ra _{\rm def} $.
\end{fac}
This result is based on Rubin's work \cite{Rub}.
We note that the structures from $\la \CC^{\rm lo} \ra _{\rm def}$ are called {\it monomorphic} by Fra\"{\i}ss\'{e}
and that $\la \CC^{\rm lo}_{\rm lab} \ra _{\rm def} $ is the class of structures admitting a {\it finite monomorphic decomposition} (FMD structures) introduced by Pouzet and Thi\'{e}ry \cite{PT}.

The next example is related to an isomorphism-closed class $\CC$ of partial orders  satisfying VC$^\sharp$
and its closure $\la \CC \ra _{\du \Pi }$ under finite products and disjoint unions.
So, if $\CC ^{\rm ba}$ is the class of Boolean algebras,
$\CC ^{\rm tree}_{\rm 0,fmd}$ the class of rooted FMD trees,
$\CC ^{\rm tree}_{\rm if,{\rm VC}^{\sharp}}$ the class of initially finite trees\footnote{A rooted tree is called initially finite iff deleting its root we obtain finitely many connectivity components.}
satisfying  {\rm VC}$^{\sharp}$,
and $\la \CC ^{\rm lo}\ra_{\dot{\cup}_{\infty}}$ is the class of infinite disjoint unions of linear orders,
then by \cite{Ksharp} we have
\begin{fac}\label{T719}
VC$^\sharp$ is true for all partial orders from the class
$$
\CC':=\la \CC ^{\rm lo}\ra _{\du\Pi }
\cup \la \CC ^{\rm ba}\ra _{\du\Pi}
\cup \la \CC ^{\rm tree}_{\rm 0,fmd}\ra _{\du\Pi }
\cup \la \CC ^{\rm tree}_{\rm if,{\rm VC}^{\sharp}}\ra _{\du\Pi }
\cup \la \CC ^{\rm lo}\ra_{\dot{\cup}_{\infty}}.
$$
\end{fac}
In this paper for a class $\CC$ of partial orders we consider its closure $\la \CC\ra _{\Sigma}$ under finite lexicographic sums.
We will say that a partial order $\X$ admits a {\it finite lexicographic decomposition with ones (largest elements)},
shortly, that $\X$ is an FLD$_1$ {\it  partial order},
iff there are a finite partial order $\BI$ and a partition $\{ X_i:i\in I\}$ of its domain $X$
such that $\X =\sum _{\BI}\X_i$ and that $\max \X _i$ exists, for each $i\in I$.
For example, each infinite linear order with a largest element has infinitely many such decompositions (into intervals of the form $(\cdot ,a]$ or $(a,b]$)
and each partial order $\X$ with a largest element has a 1-decomposition $\X =\sum _1\X$, which is trivial in our context.

In Section \ref{S3} we establish the notion of an FLD$_1$ {\it theory of partial order},
showing that a complete theory of partial order $\CT$ has an FLD$_1$ model iff all models of $\CT$ are FLD$_1$ partial orders.
By $\CD (\CT)$ we denote the class of all FLD$_1$ decompositions of models of $\CT$
and call $\sum _{\BI}\X_i\in \CD (\CT)$ a VC-{\it decomposition} (resp.\ a VC$^\sharp$-{\it decomposition})
iff $\X _i$ satisfies VC (resp.\ VC$^\sharp$), for each $i\in I$.
Then we show that VC is true for $\CT$ iff $\CT$ is large or its atomic model $\X^{\rm at}$ has a VC decomposition.
(Otherwise, $\prod _{i\in I}I(\X _i^{\rm at})=\o _1 <\c$, for each decomposition of $\X^{\rm at}$, and $I(\CT)=\o_1$; that is, $\CT$ is a counterexample.)

If $\CC$ is a class of partial orders with a largest element
for which VC is already confirmed,
in order to confirm VC for its closure $\la \CC \ra_{\Sigma}$ under finite lexicographic sums
in Section \ref{S4} we define an FLD$_1$ theory $\CT$ to be {\it actually Vaught's}
iff for some $\sum _{\BI}\X_i\in \CD (\CT)$ there are sentences $\t _i\in \Th (\X _i)$, for $i\in I$,
providing VC (e.g.\ if $\X _i$-s are reversed rooted trees).
Then we show that VC is true for each actually Vaught's FLD$_1$ theory.
Similarly, if  there exist a VC$^\sharp$-decomposition $\sum _{\BI}\X_i\in \CD (\CT)$,
then $\CT$ satisfies VC$^\sharp$.
All the aforementioned statements have natural duals
when we define FLD$_0$ partial orders (where the summands have a smallest element), FLD$_0$ theories, etc.

In Section \ref{S5} we apply these results.
First, on the basis of the results of Rubin, Steel, Schmerl and Iverson (see Fact \ref{T721}),
we confirm VC for the partial orders from the closure $\la \CC ^{\rm fin}_0 \cup \CC ^{\rm reticle}_0\cup \CC ^{\rm ba}\ra_{\Sigma}$
(their theories are actually Vaught's; note that $\CC ^{\rm lo}_0\subset\CC ^{\rm tree}_0\subset \CC ^{\rm reticle}_0$).
Second, in order to extend the result from Fact \ref{T719},
we define the operation of lexicographic sum of rooted summands, $\sum ^r$, and the corresponding closure $\la \CC\ra_{\Sigma ^r}$.
Then we show that for the class $\CC ^{{\rm VC}^{\sharp}}$ of all partial orders satisfying {\rm VC}$^{\sharp}$
we have $\la \CC ^{{\rm VC}^{\sharp}}\ra _{\Sigma ^r}=\CC ^{{\rm VC}^{\sharp}}$.
In particular {\rm VC}$^{\sharp}$ is true for each partial order from the class $\la\CC ' \ra _{\Sigma ^r}$, where
$\CC'$ is the class defined in Fact \ref{T719}.
So we obtain a large zoo of partial orders of the form
$\sum _{\BI}(\dot{\bigcup}_{j<n_i}\prod _{k<m_i^j}\X _i^{j,k})_r$
satisfying {\rm VC}$^{\sharp}$.
\section{Preliminaries}\label{S2}
\paragraph{Notation}
By $L_b$ we denote the language $\la R\ra$, where $R$ is a binary relational symbol or $\leq$, when we work with partial orders.
$\Mod _{L_b}$ denotes the class of $L_b$-structures,
and for a complete $L_b$-theory $\CT$ with infinite models
$I(\CT ,\o):=|\Mod _{L_b}(\CT,\o)/\!\cong|$ is the number of non-isomorphic countable models of $\CT$.
For simplicity, instead of $I(\CT ,\o)$ we write $I(\CT)$
and, for an $L_b$-structure $\X$, instead of $I(\Th (\X),\o)$ we write $I(\X)$;
if $\X$ is a finite structure, for convenience we define $I(\X)=1$.
For $\X \in \Mod _{L_b}$ by $\IT(\X)$ we denote the isomorphism type of $\X$ (the class of all $L_b$-structures isomorphic to $\X$)
and by $\Aut (\X )$ its automorphism group. Locally used specific notation will be explained locally.
\paragraph{Substructures on parametrically definable domains. Partitions.}
If $\f (w_0,\dots,w_{m-1},v)=\f (\bar w,v)$ is an $L_b$-formula,
then the corresponding relativization of an $L_b$-formula $\p (v_0\dots v_{n-1})=\p (\tilde v)$, where $n\in\o$,
is the $L_b$-formula $\p^{\f }(\bar w,\tilde v)$ defined by recursion in the following way
\begin{eqnarray}
\p^{\f }(\bar w,\tilde v)                                       & := & \p (\tilde v), \mbox{ if $\p (\tilde v)$ is atomic},\label{EQ727}\\
(\neg \p (\tilde v))^{\f }(\bar w,\tilde v)                     & := & \neg \p^{\f } (\bar w,\tilde v),\label{EQ728}\\
(\p_0 (\tilde v)\land \p_1 (\tilde v))^{\f }(\bar w,\tilde v)   & := & \p_0^{\f } (\bar w,\tilde v)\land \p_1^{\f } (\bar w,\tilde v),\label{EQ729}\\
(\p_0 (\tilde v)\lor \p_1 (\tilde v))^{\f }(\bar w,\tilde v)    & := & \p_0^{\f } (\bar w,\tilde v)\lor \p_1^{\f } (\bar w,\tilde v),\label{EQ730}\\
(\forall u \;\p (\tilde v,u))^{\f }(\bar w,\tilde v)            & := & \forall u\; (\f (\bar w,u)\Rightarrow \p ^{\f }(\bar w,\tilde v,u)),\label{EQ731}\\
(\exists u \;\p (\tilde v,u))^{\f }(\bar w,\tilde v)            & := & \exists u\; (\f (\bar w,u)\land \p ^{\f }(\bar w,\tilde v,u)).\label{EQ732}
\end{eqnarray}
If $\X \in \Mod _{L_b}$ and $\bar a \in X^m$, let $D_{\f (\bar a,v),\X}:=\{ x\in X: \X\models \f [\bar a,x]\}$,
let $\X \upharpoonright D_{\f (\bar a,v),\X}$ be the corresponding substructure of $\X$ and
\begin{equation}\label{EQ758}
w_{\f (\bar w,v)}(\X)  :=  \{ \IT (\X \upharpoonright D_{\f (\bar a,v),\X}) : \bar a\in X^m \}.
\end{equation}
\begin{fac}\label{T704}
If $\f (\bar w,v)$ is an $L_b$-formula,
then for each $L_b$-formula $\p (\tilde v)$, $\X \in \Mod _{L_b}$ and $\bar a \in X^m$ we have
\begin{equation}\label{EQ726}
\forall \tilde y \in (D_{\f (\bar a,v),\X})^n \;\; (\X \upharpoonright D_{\f (\bar a,v),\X}\models \p[\tilde y]  \;\;\Leftrightarrow\;\; \X \models \p^{\f }[\bar a,\tilde y]) .
\end{equation}
So, for each $L_b$-sentence $\p$, $\X \in \Mod _{L_b}$ and $\bar a \in X^m$  we have
$\X \upharpoonright D_{\f (\bar a,v),\X}\models \p  $ iff $\;\X \models \p^{\f }[\bar a]$.
\end{fac}
\dok
Permuting the universal quantifiers we fix $\X \in \Mod _{L_b}$ and $\bar a \in X^m$ and by induction show that for each $L_b$-formula $\p (\tilde v)$ we have (\ref{EQ726}).

If $\p (\tilde v)$ is an atomic formula ($v_i=v_j$ or $R(v_i,v_j)$),
$n>i,j$ and $\tilde y \in (D_{\f (\bar a,v),\X})^n$,
then, clearly, $\X \upharpoonright D_{\f (\bar a,v),\X} \models \p[y_i,y_j]$
iff $\X \models \p[y_i,y_j]$
iff (by (\ref{EQ727})) $\X \models \p ^{\f}[y_i,y_j]$;
so, (\ref{EQ726}) is true for $\p (\tilde v)$.

Assuming that  (\ref{EQ726}) is true for $\p (\tilde v)$
we prove that it is true for $\neg \p (\tilde v)$.
So, for $\tilde y \in (D_{\f (\bar a,v),\X})^n$ we have
$\X \upharpoonright D_{\f (\bar a,v),\X} \models (\neg \p(\tilde v))[\tilde y]$,
iff $\X \upharpoonright D_{\f (\bar a,v),\X} \not\models \p[\tilde y]$,
iff (by the induction hypothesis) $\X \not\models \p^{\f }[\bar a,\tilde y]$,
iff $ \X \models \neg \p^{\f }[\bar a,\tilde y]$,
iff (by (\ref{EQ728})) $ \X \models (\neg \p (\tilde v))^{\f }[\bar a,\tilde y]$.
Thus (\ref{EQ726}) is true for $\neg \p (\tilde v)$.

Assuming that  (\ref{EQ726}) is true for $\p _0(\tilde v)$ and $\p _1(\tilde v)$
we prove that it is true for $\p _0 (\tilde v)\land \p _1(\tilde v)$.
So, for $\tilde y \in (D_{\f (\bar a,v),\X})^n$ we have
$\X \upharpoonright D_{\f (\bar a,v),\X} \models (\p _0 (\tilde v)\land \p _1(\tilde v))[\tilde y]$,
iff $\X \upharpoonright D_{\f (\bar a,v),\X} \models \p_0[\tilde y]$ and $\X \upharpoonright D_{\f (\bar a,v),\X} \models \p_1[\tilde y]$,
iff (by the induction hypothesis) $\X \models \p_0^{\f }[\bar a,\tilde y]$ and $\X \models \p_1^{\f }[\bar a,\tilde y]$,
iff $\X \models (\p_0^{\f }\land \p_1^{\f })[\bar a,\tilde y]$,
iff (by (\ref{EQ729})) $ \X \models (\p _0 (\tilde v)\land \p _1(\tilde v))^{\f }[\bar a,\tilde y]$.
Thus (\ref{EQ726}) is true for $\p _0 (\tilde v)\land \p _1(\tilde v)$
and, similarly, for $\p _0 (\tilde v)\lor \p _1(\tilde v)$.

Assuming that  (\ref{EQ726}) is true for $\p (\tilde v,u)$
we prove that it is true for $\forall u\;\p (\tilde v,u)$.
Let $Y:=D_{\f (\bar a,v),\X}$ and $\tilde y \in Y^n$; then
$\X \upharpoonright Y \models (\forall u\;\p (\tilde v,u))[\tilde y]$,
iff for each $y\in Y$ we have $\X \upharpoonright Y \models \p (\tilde v,u)[\tilde y,y]$,
iff (by the induction hypothesis) for each $y\in Y$ we have $\X \models \p ^{\f }(\bar w,\tilde v,u)[\bar a,\tilde y,y]$
iff for each $y\in X$ we have that $\X\models \f [\bar a,y]$ implies $\X \models \p ^{\f }(\bar w,\tilde v,u)[\bar a,\tilde y,y]$
iff for each  $y\in X$ we have $\X\models (\f (\bar w,u)\Rightarrow \p ^{\f }(\bar w,\tilde v,u))[\bar a,\tilde y,y]$
iff $\X\models \forall u \;(\f (\bar w,u)\Rightarrow \p ^{\f }(\bar w,\tilde v,u))[\bar a,\tilde y]$
iff (by (\ref{EQ731})) $\X\models  (\forall u \;\p (\tilde v,u))^{\f }  [\bar a,\tilde y]$.
So, (\ref{EQ726}) is true for $\forall u\;\p (\tilde v,u)$ and, similarly, for $\exists u\;\p (\tilde v,u)$.
\kdok
\begin{fac}\label{T706}
If $\X$ is an $L _b$-structure, $\f (\bar w,v)$ an $L_b$-formula, $\bar a\in X^m$, where $D_{\f (\bar a,v),\X}\neq \emptyset$ and $\X \preccurlyeq \Y$,
then
\begin{equation}\label{EQ759}
D_{\f (\bar a, v),\X}=X\cap  D_{\f (\bar a, v),\Y} \;\;\mbox{ and }\;\; \X \upharpoonright D_{\f (\bar a, v),\X}\preccurlyeq \Y \upharpoonright D_{\f (\bar a, v),\Y}.
\end{equation}
\end{fac}
\dok
First, $x\in D_{\f (\bar a, v),\X}$
iff $x\in X$ and $\X\models \f [\bar a, x]$
iff (since $\X \preccurlyeq \Y$) $x\in X$ and $\Y\models \f [\bar a, x]$
iff $x\in X \cap D_{\f (\bar a, v),\Y}$;
so, $ D_{\f (\bar a, v),\X}=X\cap  D_{\f (\bar a, v),\Y}$.
Second, for a formula $\p (\tilde w)$ and $\tilde x \in (D_{\f (\bar a, v),\X})^n$ we have
$\X \upharpoonright D_{\f (\bar a, v),\X} \models \p [\tilde x]$,
iff (by Fact \ref{T704}) $\X \models \p ^{\f }[\bar a,\tilde x]$,
iff (since  $\X \preccurlyeq \Y$) $\Y \models \p ^{\f }[\bar a,\tilde x]$,
iff (by Fact \ref{T704}) $\Y \upharpoonright D_{\f (\bar a, v),\Y} \models \p [\tilde x]$.
Thus, $\X \upharpoonright D_{\f (\bar a, v),\X}\preccurlyeq \Y \upharpoonright D_{\f (\bar a, v),\Y}$.
\kdok
\begin{fac}\label{T701}
If $\X \cong \Y$, then $w_{\f (\bar w,v)}(\X)=w_{\f (\bar w,v)}(\Y)$. If $|X|=\o$, then $w_{\f (\bar w,v)}(\X)\leq \o$.
\end{fac}
\dok
Let $f:\X \rightarrow \Y$ be an isomorphism.
If $\t \in w_{\f (\bar w,v)}(\X)$,
then by (\ref{EQ758}) there is $\bar a\in X^m$ such that $\t =\IT(\X \upharpoonright D_{\f (\bar a,v),\X})$.
For $x\in X$ we have $x\in D_{\f (\bar a,v),\X}$
iff $\X \models \f [\bar a,x]$
iff $\Y \models \f [f\bar a,fx]$
iff $f(x)\in D_{\f (f\bar a,v),\Y}$.
Thus $f[D_{\f (\bar a,v),\X}]=D_{\f (f\bar a,v),\Y}$
and, hence, $\X \upharpoonright D_{\f (\bar a,v),\X} \cong \Y \upharpoonright D_{\f (f\bar a,v),\Y}$,
that is $\t=\IT(\X \upharpoonright D_{\f (\bar a,v),\X}) =\IT ( \Y \upharpoonright D_{\f (f\bar a,v),\Y})\in w_{\f (\bar w,v)}(\Y)$.
So, $w_{\f (\bar w,v)}(\X)\subset w_{\f (\bar w,v)}(\Y)$
and the proof of the other inclusion is symmetric.
If $|X|=\o$, then $|X^m|=\o$ and, by (\ref{EQ758}),  $w_{\f (\bar w,v)}(\X)\leq \o$.
\kdok
Concerning partitions of structures by Proposition 2.3 of \cite{Ksharp} we have
\begin{fac}\label{T216}
If $\X$ is a countable $L_b$-structure and $\{ X_i : i\in I\}$ a partition of its domain $X$, then

(a) If for each $f\in \Aut (\X  )$ and $i\in I$ from $f[X_i]\cap X_i \neq \emptyset$ it follows that $f[X_i]=X_i$,
then
\begin{equation} \label{EQ237}\textstyle
\X  \mbox{ is $\o$-categorical} \Rightarrow \forall i\in I \;(\X_i \mbox{ is $\o$-categorical});
\end{equation}

(b) If $|I|<\o$ and  $\bigcup _{i\in I}f_i \in \Aut (\X )$, whenever $f_i \in \Aut (\X_i )$, for $i\in I$,
then we have ``$\Leftarrow$" in (\ref{EQ237}).
\end{fac}
\paragraph{Lexicographic sums of $L_b$-structures}
Let $\BI=\la  I, \r _I\ra$ and $\X _i =\la X_i , \r _i \ra$, $i\in I$, be $L_b$-structures with pairwise disjoint domains.
The {\it lexicographic sum of the structures $\X _i$, $i\in I$, over the structure $\BI$}, in notation $\sum _{\BI}\X _i$, is the  $L_b$-structure
$\X :=\la X, \r \ra$, where $X := \bigcup _{i\in I}X_i$ and
for $x,x'\in X$ we have: $x\,\r \,x'$ iff
\begin{equation}\label{EQ200}
\exists i\in I \;\; \Big(x,x'\in X_i \land x \,\r _i\, x'\Big) \;\lor\;
\exists \la i,j \ra \in \r _I \setminus \Delta _I \;\; \Big(x\in X_i \land x'\in X_j \Big),
\end{equation}
where $\Delta _I:=\{ \la i,i\ra:i\in I\}$.
If $\X$ and $\Y$ are $L_b$-structures, by $\PI(\X ,\Y )$ we denote the set of all partial isomorphisms between $\X$ and $\Y$.
Let $\EF _k (\X ,\Y )$ denote that Player II has a winning strategy in the Ehrenfeucht-Fra\"{\i}ss\'{e} game of length $k$ between $\X $ and $\Y $.
\begin{fac}\label{T200}
Let $\sum _{\BI}\X _i$ and $\sum _{\BI}\Y _i$ be lexicographic sums of $L_b$-structures.
Then

(a) If $f_i \in \PI(\X _i , \Y _i)$, for $i\in J\subset I$, then $\bigcup _{i\in J}f_i \in \PI(\sum _{\BI}\X _i ,  \sum _{\BI}\Y _i) $;

(b) If $\X _i \preccurlyeq  \Y _i$, for all $i\in I$, then $\sum _{\BI}\X _i \preccurlyeq \sum _{\BI}\Y _i $;

(c) If $\X _i \equiv  \Y _i$, for all $i\in I$, then $\sum _{\BI}\X _i \equiv  \sum _{\BI}\Y _i $;

(d) If $|I|<\o$ and $I(\X _i)=1$, for all $i\in I$, then $I(\sum _{\BI}\X _i )=1$.
\end{fac}
\dok
(a) This claim follows directly from (\ref{EQ200}).

(b) First we recall a standard fact:
If $\X$ and $\Y$ are $L$-structures, where $|L|<\o$, and $\X \subset \Y$, then
\begin{equation}\label{EQ533}
\X \preccurlyeq \Y \;\;\mbox{ iff }\;\;\forall n\in \o \;\;\forall \bar x \in  X^n\;\;\forall k\in \o \;\; \EF _k ((\X ,\bar x),(\Y ,\bar x)) .
\end{equation}
Namely, if $\X \subset \Y$,
then (see \cite{Bell}, p.\ 77) $\X \preccurlyeq \Y$
iff for each $n\in \o$ and $\bar x \in  X^n$ we have $(\X ,\bar x)\equiv  (\Y ,\bar x)$,
where $\equiv$ refers to the language $L_{\bar c}=L \cup \{ c_0, \dots ,c_{n-1}\}$ and $c_i$, $i<n$, are new constants.
Since $|L_{\bar c}|<\o$, we have $(\X ,\bar x)\equiv (\Y ,\bar x)$ iff $\EF _k ((\X ,\bar x),(\Y ,\bar x))$, for all $k\in \o$.

Now, if $\X _i \preccurlyeq  \Y _i$, for $i\in I$,
then for each $i\in I$ we have
\begin{equation}\label{EQ582}
\forall n\in \o \;\;\forall \bar x \in  X_i^n\;\;\forall k\in \o \;\; \EF _k ((\X _i ,\bar x),(\Y _i ,\bar x))
\end{equation}
and we have to prove that
for each $n\in \o$, $\bar x =\la x_j :j<n\ra \in  (\bigcup _{i\in I}X_i)^n$ and $k\in \o$ we have $\EF _k ((\sum _{\BI}\X _i ,\bar x),(\sum _{\BI}\Y _i ,\bar x))$.
For $i\in I$ let $J_i:=\{ j<n : x_j \in X_i \}$;
then $\bar x ^i:=\la x_j :j\in J_i\ra \in X_i ^{|J_i|}$ (possibly $J_i= \bar x ^i=\emptyset$)
and by (\ref{EQ582}) Player II has a winning strategy in the game $\EF _k ((\X _i ,\bar x ^i),(\Y _i ,\bar x ^i))$, say $\Sigma _i$.
In the game $\EF _k ((\sum _{\BI}\X _i ,\bar x),(\sum _{\BI}\Y _i ,\bar x))$ Player II uses the strategies $\Sigma _i$, $i\in I$,
in a natural way: when Player I chooses an element $a$ from $X_i$ or $Y_i$
Player II responds applying $\Sigma _i$
to the restriction of the previous play between $(\X _i ,\bar x ^i)$ and $(\Y _i ,\bar x ^i)$ extended by $a$.
In this way
at the end of the game for each $i\in I$ a partial isomorphism $f_i\in \PI ((\X _i ,\bar x ^i),(\Y _i ,\bar x ^i))$ is obtained
(some of them can be empty)
by (a) we have $f:=\bigcup _{i\in I}f_i \in \PI(\sum _{\BI}\X _i ,  \sum _{\BI}\Y _i) $.
Moreover, by the construction we have $\{ \la x_j ,x_j \ra : j<n\}\cup f \in \PI ((\sum _{\BI}\X _i ,\bar x),(\sum _{\BI}\Y _i ,\bar x))$,
and Player II wins.

(c) Let $\X _i \equiv  \Y _i$, for all $i\in I$.
For $k\in\N$ and $i\in I$, let $\Sigma _i$ be a winning strategy for Player II in the game $\EF _k (\X _i , \Y _i)$.
If in the game $\EF _k (\sum _{\BI}\X _i ,\sum _{\BI}\Y _i)$ Player II follows $\Sigma _i$'s
(when Player I chooses an element from $X_i$ or $Y_i$, Player II uses $\Sigma _i$),
then, by (a), Player II wins.
So $\sum _{\BI}\X _i \equiv _k\sum _{\BI}\Y _i$ for all $k\in\o$ and, hence, $\sum _{\BI}\X _i \equiv  \sum _{\BI}\Y _i $.

(d) Let $|I|<\o$ and let $I(\X _i)=1$, for all $i\in I$.
If $f_i\in \Aut (\X _i)$, for $i\in I$,
then by (a) we have $\bigcup _{i\in I}f_i\in \Aut (\X)$.
By Fact \ref{T216}(b) we have $I(\X )=1$.
\kdok
\section{FLD theories}\label{S3}
We will say that a partial order $\X$ admits a {\it finite lexicographic decomposition with ones (largest elements)},
shortly, that $\X$ is an FLD$_1$-{\it partial order},
iff there are a finite partial order $\BI=\la I,\leq _{\BI}\ra$ and a partition $\{ X_i:i\in I\}$ of its domain $X$
such that $\X =\sum _{\BI}\X_i$ and that $\max \X _i$ exists, for each $i\in I$.
\begin{te}\label{T700}
If $\X$ is an FLD$_1$-partial order
and $\X =\sum _{\BI}\X_i$,
where $\BI=\la n,\leq _{\BI}\ra$ and $\X_i=\la X_i,\leq _i\ra$, for $i<n$, are pairwise disjoint partial orders,
$r_i=\max \X _i$, for $i<n$, and $\bar r:=\la r_0,\dots ,r_{n-1}\ra$, then \\[-8mm]
\begin{itemize}\itemsep=-1mm
\item[(a)] $X_i=D_{\f _i (\bar r , v), \X}$, for $i<n$, where
               \begin{equation}\label{EQ706}\textstyle
               \f _i(\bar r , v):= v\leq r_i \land \bigwedge _{j\,<_\BI\, i}r_j <v \land \bigwedge _{j\,>_\BI\, i}r_j >v \land \bigwedge _{j\;\not\parallel _{\BI} \; i}r_j \not\parallel v;
               \end{equation}
\item[(b)] The formula $\ve (\bar r, u,v)=\bigvee _{i<n} (\f _i(\bar r , u)\land \f _i(\bar r , v))$
               defines in $\X$ an equivalence relation on the set $X$ and $X/D_{\ve (\bar r, u,v),\X}=\{ X_i :i<n\}$;

\item[(c)] If  $\,\Y  \equiv \X$ and $\t _i \in \Th (\X _i)$, for $i<n$,\footnote{Or, equivalently, if $\CT _i \in [\Th (\X _i)]^{<\o}$ and $\t _i=\bigwedge \CT _i$, for $i<n$.}
               then there is $\bar r ':=\la r_0',\dots ,r_{n-1}'\ra \in Y^n$ such that defining $Y_i:=D_{\f _i (\bar r ' , v), \Y}$, for $i<n$, we have\\[-8mm]
               \begin{itemize}\itemsep=-1mm
               \item[(i)] $\{ Y_i:i<n\}$ is a partition of the set $Y$  and $Y/D_{\ve (\bar r',u,v),\Y} =\{ Y_i:i<n\}$,
               \item[(ii)] $\Y =\sum _{\BI}\Y_i$ and $r_i'=\max \Y _i$, for $i<n$,
               \item[(iii)] $\Y _i \models \t _i$, for $i<n$.
               \end{itemize}
\end{itemize}
\end{te}
\dok
(a) Let $\X =\la X, \leq\ra$. The sets $X_i$, $i<n$, are pairwise disjoint and by (\ref{EQ200}) for  $x,y\in X$ we have
\begin{equation}\label{EQ707}
x\leq y \Leftrightarrow
\exists i\in n \;\; (x,y\in X_i \land x \,\leq _i\, y) \;\lor\;
\exists  i,j \in n \;\; ( i<_\BI j \land x\in X_i \land y\in X_j ).
\end{equation}
We take $i<n$ and show that for each $x\in X$ we have $x \in X_i$ iff $\X \models \f_i [\bar r , x]$, namely,
\begin{equation}\label{EQ708}\textstyle
x \in X_i \;\Leftrightarrow\;
x\leq r_i \land \bigwedge _{j<_\BI i}r_j <x \land \bigwedge _{j>_\BI i}r_j >x \land \bigwedge _{j\;\not\parallel _\BI \; i}r_j \not\parallel x.
\end{equation}
If $x\in X_i$,
then $x,r_i=\max \X _i\in X_i$
and, hence, $x \leq _i r_i$,
which by (\ref{EQ707}) gives $x \leq  r_i$.
If $j<_\BI i$,
then, since $r_j\in X_j$ and $x\in X_i$,
by (\ref{EQ707}) we have $r_j <x$.
If $j>_\BI i$,
then, since $r_j\in X_j$ and $x\in X_i$,
by (\ref{EQ707}) we have $r_j >x$.
Finally, if $j\not\parallel _\BI  i$,
then $i\neq j$ and $x \neq r_j$.
Assuming that $x<r_j$
by (\ref{EQ707}) we would have $i<_\BI j$,
which is false.
So, $x\not<r_j$
and, similarly, $r_j\not< x$,
which gives $r_j \not\parallel x$. So, ``$\Rightarrow$" in (\ref{EQ708}) is proved.

Let the r.h.s.\ of (\ref{EQ708}) be true; we prove that $x\in X_i$.
Assuming that $x\in X_j$, for some $j\neq i$, we would have $x\leq _j r_j=\max \X _j$
and by (\ref{EQ707}) $x\leq r_j$.
Clearly we have $i<_\BI j$ or $j<_\BI i$ or $i\not\parallel _\BI j$.
Now, if $i<_\BI j$,
then, since $r_i \in X_i$,
by (\ref{EQ707}) we would have $r_i < x$,
which is false because by the r.h.s.\ of (\ref{EQ708}) $x\leq r_i$.
If $j<_\BI i$,
then by the r.h.s.\ of (\ref{EQ708}) we would have $r_j <x$,
which is false because $x\leq r_j$.
Finally, if $i\not\parallel _\BI j$
then by the r.h.s.\ of (\ref{EQ708}) we would have $r_j \not\parallel x$
which is false  because $x\leq r_j$.
Thus $x\in X_i$ and (\ref{EQ708}) is proved.
So, $X_i=D_{\f _i (\bar r , v), \X}$, for $i<n$.

(b) Since $\{ X_i :i<n\}$ is a partition of the set $X$ claim (b) follows from (a).

(c) By (a) we have $X=\bcd _{i<n}D_{\f _i (\bar r , v), \X}$.
So $\X \models \f _p [\bar r)]$,
where $\f _p (\bar w)$ is the $L_b$-formula saying that $\{ D_{\f _i (\bar w,v), \cdot}: i<n\}$ is a partition of the domain; say,
\begin{equation}\label{EQ714}\textstyle
\f _p (\bar w):= \forall v \;\;\bigvee _{i<n} (\f _i (\bar w,v)\land \bigwedge _{j\in n\setminus \{i\}}\neg \f _j (\bar w,v)).
\end{equation}
Since $\X =\sum _{\BI}\X _i$ by (\ref{EQ200}) we have: if $i,j<n$, $i\neq j$, $x\in X_i$ and $x'\in X_j$, then $x \leq _\X x'$ iff $i<_\BI j$.
Thus, by (a), $\X\models \f _m [\bar r]$,  where
\begin{eqnarray}\textstyle
\f _m (\bar w) & := & \textstyle\bigwedge _{i,j\in n \, \land \, i<_\BI j}\forall u,v \;(\f _i(\bar w ,u)\land \f _j(\bar w ,v)\Rightarrow u \leq v) \;\land\;\nonumber \\
               &    &  \textstyle\bigwedge _{i,j\in n \, \land \, i\not\leq_\BI j}\forall u,v \; (\f _i(\bar w ,u)\land \f _j(\bar w ,v)\Rightarrow \neg u \leq v).\label{EQ716}
\end{eqnarray}
In $\X$ for $i<n$ we have $r_i=\max D_{\f _i (\bar r , v), \X}$;
so, for each $x\in D_{\f _i (\bar r , v), \X}$ we have $x\leq _\X r_i $,
that is, $\X\models \f _{r_i}[\bar r]$,
where
\begin{equation}\label{EQ721}
\f _{r_i}(\bar w):= \forall v \;( \f _i (\bar w , v)\Rightarrow v\leq  w_i).
\end{equation}
Let $i<n$ and $\t _i \in \Th (\X _i)$.
Then by (a) $\X _i =\X \upharpoonright D_{\f _i (\bar r , v), \X}\models \t_i$
and by Fact \ref{T704} $\X \models \t_i ^{\f _i } [\bar r]$.
Moreover, $\X \models \f _{\CT }[\bar r]$, where
\begin{equation}\label{EQ720}\textstyle
\f _{\CT }(\bar w):=\bigwedge _{i<n}\t_i ^{\f _i } (\bar w).
\end{equation}
Thus we have $\X \models \f [\bar r]$,
where $\f (\bar w):= \f _p(\bar w)\land \f _m(\bar w)\land \bigwedge _{i<n}\f _{r_i}(\bar w)\land\f _{\CT }(\bar w)$
and, hence, $\X\models \exists \bar w  \;\f (\bar w)$.
Consequently, $\Y\models \exists \bar w \; \f (\bar w)$
so there is $\bar r '\in Y^n$ such that $\Y\models \f [\bar r']$
and we check (i)--(iii).

(i) Since $\Y\models \f _p [\bar r']$ and $Y_i:=D_{\f _i (\bar r',v), \Y}$, for $i<n$, by (\ref{EQ714}) $\{ Y_i: i<n\}$ is a partition of $\Y$.

(ii) Since $\Y\models \f _m [\bar r']$ and $Y_i:=D_{\f _i (\bar r',v), \Y}$, for $i<n$, by (\ref{EQ716}) we have
\begin{equation}\label{EQ717}\textstyle
\forall i,j<n \;\;(i\neq j \Rightarrow \forall y\in Y _i \;\;\forall y'\in Y_j \;\;(y\leq _\Y y' \Leftrightarrow i<_\BI j)).
\end{equation}
Since $\X$ is a partial order and $\Y  \equiv \X$,  $\Y =\la Y,\leq _\Y\ra$ is a partial order;
so, for each $i<n$ its substructure $\Y _i=\la Y_i ,\leq _{\Y _i}\ra$, where $\leq _{\Y _i}:=\leq _\Y \cap Y_i^2$, is a partial order.
For a proof that $\Y =\sum _{\BI}\Y_i $ we have to show that for each $y,y'\in Y$ we have
\begin{equation}\label{EQ718}
y \leq _\Y y'  \Leftrightarrow
\exists i<n \;\; (y,y'\in Y_i \land y\leq _{\Y _i} y') \;\lor\;
\exists  i,j <n \;\; (y\in Y_i \land y'\in Y_j \land i<_\BI j ).
\end{equation}
Let $y \leq _\Y y'$.
If $y,y'\in Y_i $, for some $i<n$,
then, since $\Y_i $ is a substructure of $\Y$,
we have $y\leq _{\Y _i} y'$ and the r.h.s.\ of (\ref{EQ718}) is true.
Otherwise, there are different $i,j <n$ such that $y\in Y_i$ and $y'\in Y_j$;
so, by (\ref{EQ717}) we have $i<_\BI j$ and the r.h.s.\ of (\ref{EQ718}) is true again.
Conversely, let the r.h.s.\ of (\ref{EQ718}) be true.
If $y,y'\in Y_i$ and $y\leq _{\Y _i} y'$,
then, since $\Y _i \subseteq \Y$, we have $y \leq _\Y y'$.
Otherwise we have $y\in Y_i$ and $y'\in Y_j$, where $i<_\BI j$
and, by (\ref{EQ717}), $y \leq _\Y y'$ again.
So, (\ref{EQ718}) is true and $\Y =\sum _{\BI}\Y_i $.
In addition,
for each $i<n$ we have $\Y\models \f _{r_i}[\bar r ']$
and, since $Y_i:=D_{\f _i (\bar r ' , v), \Y}$,
by (\ref{EQ721}) we have $r_i'=\max \Y _i$.

(iii) Since $\Y\models \f _{\CT} [\bar r']$,
by (\ref{EQ720}) for $i<n$
we have $\Y \models \t_i^{\f _i } [\bar r ']$
and by Fact \ref{T704} $\Y _i \models \t_i$.
\kdok
For an FLD$_1$-partial order $\X$ let
$\I (\X):=\{ \BI \in \CC ^{\rm fin}: \mbox{ there is an FLD$_1$ decomposition } \X =\sum _{\BI}\X_i\}$,
where $\CC ^{\rm fin}$ is the class of finite partial orders.
Clearly each poset $\BI \in \CC ^{\rm fin}$ is isomorphic to one with domain $n$, for some $n\in \N$;
so, when it is convenient we can assume that $\BI=\la n,\leq _{\BI}\ra$.
\begin{te}\label{T708}
If $\CT$ is a complete theory of partial order, then

(a) If $\X$ is an FLD$_1$-model of $\CT$, then $\I (\Y)=\I (\X)$, for each model $\Y$ of $\CT$;

(b) $\CT$ has an FLD$_1$-model iff all models of $\CT$ are FLD$_1$-partial orders.
\end{te}
\dok
If $\X \models \CT$ is an FLD$_1$-partial order, $\BI \in \I (\X)$ and $\Y \equiv \X$,
then by Theorem \ref{T700}(c) there is an $\BI$ decomposition $\Y =\sum _{\BI}\Y _i$ of $\Y$;
thus $\BI \in \I (\Y)$ and $\Y$ is an FLD$_1$-partial order.
So, $\I (\X)\subset\I (\Y)$ and, analogously, $\I (\Y)\subset\I (\X)$.
Thus (a) is true and (b) follows from (a).
\kdok
According to Theorem \ref{T708}(b)  a complete theory of partial order $\CT$ will be called an FLD$_1$-theory
iff some model of $\CT$ is an FLD$_1$-partial order (iff all models of $\CT$ are FLD$_1$-p.o.-s).
Then, by Theorem \ref{T708}(a), we legally define $\I(\CT):= \I (\X)$, where $\X$ is some (any) model of $\CT$.

In addition, if $\X$ is an FLD$_1$-partial order, $\BI\in \I (\X)$ and $\X =\sum _{\BI}\X_i$ is an FLD$_1$ decomposition of $\X$,
then it is possible that there are more such $\BI$-decompositions $\X =\sum _{\BI}\X_i'$, where $\X _i'\not\cong \X_i$ for some $i$-s (e.g.\ for linear orders).
So let $\CD _{\BI}(\X)$ be the class of all $\BI$-decompositions of $\X$,
$\CD (\X):=\bigcup _{\BI\in \I (\X)}\CD _{\BI}(\X)$ the class of all FLD$_1$-decompositions of $\X$
and $\CD (\Th (\X))=\bigcup _{\Y \equiv \X}\CD (\Y)$ the class of all FLD$_1$-decompositions of models of $\Th (\X)$.
So, for an FLD$_1$-theory $\CT$ we define the class
\begin{equation}\textstyle\label{EQ740}
\CD (\CT):=\bigcup _{\X \models \CT}\bigcup _{\BI\in \I (\CT)}\CD _{\BI}(\X).
\end{equation}
\begin{te}\label{T710}
If $\CT$ is an FLD$_1$-theory, then the following conditions are equivalent

(a) $I(\CT)=1$, that is, $\CT$ is $\o$-categorical;

(b) For each $\sum _{\BI}\X_i\in \CD (\CT)$ we have $I(\X _i)=1$, for all $i\in I$;

(c) There is $\sum _{\BI}\X_i\in \CD (\CT)$ such that $I(\X _i)=1$, for all $i\in I$.
\end{te}
\dok
Let, in addition, $\CT$ be a theory with infinite models (otherwise, the statement is obviously true).

(a) $\Rightarrow$ (b).
Let $I(\CT)=1$,
let $\sum _{\BI}\X_i\in \CD (\CT)$, where $\X:=\sum _{\BI}\X_i$ is countable,
let $r_i:=\max \X _i$, for $i\in I$, and $\bar r:=\la r_i :i\in I\ra$.
We fix $i_0\in I$ and prove that $I(\X _{i_0})=1$; if $|X_{i_0}|<\o$ we are done; so let $|X_{i_0}|=\o$.
Since $\bar r\in X^I$ and $\X$ is $\o$-categorical,
the expansion $(\X ,\bar r)$ of $\X$ to $L_{\bar c}:=\la \leq ,\la c_i:i\in I\ra\ra$ is $\o$-categorical (see \cite{Hodg}, p.\ 346).
Defining $Y:=X_{i_0}\cup \{ r_i :i\in I\}$ and $\Y :=\la Y, \leq ^\X \upharpoonright Y\ra$
we obtain a substructure $(\Y ,\bar r)$ of $(\X ,\bar r)$.
Since $\Y =\sum _{\BI}\Y_i$, where $\Y_{i_0}=\X_{i_0}$ and $Y_i=\{ r_i\}$, for $i\neq i_0$,
by Theorem \ref{T700}(a) we have $X_{i_0}=D_{\f _{i_0}(\bar r, v), \X }$;
so, for the $L_{\bar c}$-formula $\p _{i_0}(v):=\f _{i_0}(\bar c, v)\lor \bigvee _{i\in I\setminus \{ i_0\}}v=c_i $
we have $Y=D_{\p _{i_0}(v), (\X ,\bar r)}$.
Thus the $L_{\bar c}$-structure $(\Y ,\bar r) =(\X ,\bar r) \upharpoonright D_{\p _{i_0}(v), (\X ,\bar r)}$
is the relativization of the $L_{\bar c}$-structure $(\X ,\bar r)$ to the $\emptyset$-definable set $D_{\p _{i_0}(v), (\X ,\bar r)}$,
and $\Y :=\la X_{i_0}\cup \{ r_i :i\in I\}, \leq ^\X \upharpoonright Y\ra$ is the corresponding relativized reduct of $(\X ,\bar r)$ to $L_b=\la \leq\ra$.
So, since the structure $(\X ,\bar r)$ is $\o$-categorical, $\Y$ is $\o$-categorical too (see \cite{Hodg}, p.\ 346).

In order to prove that the partial order $\X _{i_0}$ is $\o$-categorical
we apply Fact \ref{T216}(a) to $\Y$ and its partition $\{ Y_i:i\in I\}=\{ X_{i_0}\}\cup \{ \{r_i\} :i\in I\setminus \{i_0\}\}$.
So, we have to prove that for each $f\in \Aut (\Y)$ and each $i\in I$
from $f[Y_i]\cap Y_i \neq \emptyset$ it follows that $f[Y_i]=Y_i$.
First, if $i\in I\setminus \{i_0\}$ and $f[\{r_i\}]\cap \{r_i\} \neq \emptyset$,
then $f(r_i)=r_i$ and we are done;
so $i_0$ remains to be considered.
Since $|X_{i_0}|=\o$ and $|Y\setminus X_{i_0}|<\o$,
we will always have $f[X_{i_0}]\cap X_{i_0} \neq \emptyset$;
thus we have to prove that $f[X_{i_0}]= X_{i_0}$, for each $f\in \Aut (\Y)$.
First we show that
\begin{equation}\label{EQ744}\textstyle
f(r_{i_0})=r_{i_0}.
\end{equation}
Assuming that $f(r_{i_0})\neq r_{i_0}$ we have three cases (recall that $\Y =\sum _{\BI}\Y_i$).
First, if $f(r_{i_0})< r_{i_0}$,
then, since $f^{-1}\in \Aut (\Y)$ too, we would have $r_{i_0}< f^{-1}(r_{i_0})< f^{-1}(f^{-1}(r_{i_0}))<\dots$
and, hence, $|(r_{i_0},\cdot)_\Y|=\o$, which is false because $X_{i_0}\subset (\cdot ,r_{i_0}]_\Y$.
Second, if $f(r_{i_0})> r_{i_0}$, we would have $r_{i_0}< f(r_{i_0})< f(f(r_{i_0}))<\dots$,
which is false for the same reason.
Third, if $f(r_{i_0}) \not\parallel r_{i_0}$,
then, since $X_{i_0}\subset (\cdot ,r_{i_0}]_\Y$, we would have $f(r_{i_0}) = r_i$, for some $i\neq i_0$,
and, hence $r_i \not\parallel r_{i_0}$.
So, since $\Y =\sum _{\BI}\Y_i$  we would have $x \not\parallel r_i$, for all $x\in X_{i_0}$.
But taking $x\in f[X_{i_0}]\cap X_{i_0}$, since $f[X_{i_0}]\subset (\cdot ,r_i]_\Y$ we would have $x\leq r_i$,
which gives a contradiction.
Thus (\ref{EQ744}) is true.

Next we prove that $f[\{ r_i:i\neq i_0\}]=\{ r_i:i\neq i_0\}$,
which will imply that $f[X_{i_0}]= X_{i_0}$.
On the contrary, suppose that there is $i\neq i_0$ such that $f(r_i)\in X_{i_0}$;
so, by (\ref{EQ744}), $f(r_i)\in X_{i_0}\setminus \{ r_{i_0}\}$.
We have three cases again.
First, if $r_i <r_{i_0}$,
then, since $\Y =\sum _{\BI}\Y_i$ and, consequently, $r_i <X_{i_0}$,
we would have $r_i <f(r_i)$ and, hence, $r_i >f^{-1}(r_i)>f^{-1}(f^{-1}(r_i))>\dots$
and, hence, $|(\cdot ,r_i)_\Y |=\o$, which is false because $X_{i_0}\subset (r_i ,\cdot )_\Y$.
Second, if $r_i >r_{i_0}$,
we would have $f(r_i)<r_i$
and, hence, $r_i <f^{-1}(r_i)<f^{-1}(f^{-1}(r_i))<\dots$;
thus $|(r_i,\cdot )_\Y |=\o$,
which is false because $X_{i_0}\subset (\cdot ,r_i )_\Y$.
Third, if $r_i \not\parallel r_{i_0}$,
then we would have $f(r_i) \not\parallel f(r_{i_0})=r_{i_0}$,
which is false because $f(r_i)\in X_{i_0}$
and, hence, $f(r_i) \leq r_{i_0}$.
Thus $f[X_{i_0}]= X_{i_0}$ and by Fact \ref{T216}(a) the partial order $\X _{i_0}$ is $\o$-categorical.

Clearly,  (b) implies (c) and, by Fact \ref{T200}(d),  (c) implies (a).
\kdok
\begin{te}\label{T711}
If $\CT$ is an FLD$_1$-theory, $\sum _{\BI}\X_i\in \CD (\CT)$ and $I(\X _i)>\o$, for some $i\in I$, then
$$\textstyle
I(\CT)\geq \prod _{i\in I}I(\X _i)\in \{ \o _1,\c\}.
$$
\end{te}
\dok
Let $\X :=\sum _{\BI}\X_i\in \CD (\CT)$, where $\BI=\la n,\leq _{\BI}\ra$ and let $I(\X _{i_0})=\k >\o$.
Let $r_i=\max \X _i$, for $i<n$, and $\bar r:=\la r_0,\dots ,r_{n-1}\ra$.
Then $|\X|\geq \o$
and by the L\"{o}wenheim-Skolem theorem and Fact \ref{T200}(b) we can assume that $|X_i|\leq \o$, for all $i<n$.
Let $\Z ^\a$, $\a <\k$, be non-isomorphic countable models of $\Th (\X _{i_0})$
and w.l.o.g.\ assume that for each $\a <\k$ we have $\max (\Z ^\a)=r_{i_0}$
and that $Z^\a \cap X_i=\emptyset$, for all $i<n$.
For $\a <\k$, let $\Y^\a :=\sum _{\BI}\Y _i^\a$, where for $i<n$ we define
\begin{equation}\label{EQ712}\textstyle
\Y _i^\a := \left\{
           \begin{array}{cl}
           \Z ^\a, & \mbox{ if } i =i_0, \\
           \X_i,          & \mbox{ if } i\neq i_0.
           \end{array}
           \right.
\end{equation}
Then by Fact \ref{T200}(c) we have $\Y ^\a \in \Mod (\CT,\o)$,
by Theorem \ref{T700}(a) and (\ref{EQ712}) we have $Z^\a=\Y _{i_0}^\a=D_{\f _{i_0} (\bar r , v), \Y ^\a}$
and, hence, $\Z^\a= \Y ^\a \upharpoonright D_{\f _{i_0} (\bar r , v), \Y ^\a}$.

Thus $\IT(\Z^\a) \in w_{\f _{i_0}(\bar w,v)}(\Y ^\a)  :=  \{ \IT (\Y ^\a \upharpoonright D_{\f_{i_0} (\bar a,v),\Y ^\a}) : \bar a\in (Y^\a)^n \}$, for each $\a <\k$,
and, hence,
\begin{equation}\label{EQ713}\textstyle
\{ \IT(\Z^\a):\a <\k\}\subset \bigcup _{\a <\k} w_{\f _{i_0}(\bar w,v)}(\Y ^\a).
\end{equation}
For each $\a<\k$ we have $|Y^\a|=\o$
and, hence, $|w_{\f _{i_0}(\bar w,v)}(\Y ^\a)|\leq \o$.
So, since $|\{ \IT(\Z^\a):\a <\k\}|=\k >\o$,
by (\ref{EQ713}) there are $\k$-many different sets $w_{\f _{i_0}(\bar w,v)}(\Y ^\a)$,
say $w_{\f _{i_0}(\bar w,v)}(\Y ^{\a_\xi})$, $\xi<\k$,
which by Fact \ref{T701} implies that $\Y ^{\a_\xi}\not\cong\Y ^{\a_\zeta}$, for $\xi\neq \zeta$.
Thus $\Y ^{\a_\xi}$, $\xi <\k$, are non-isomorphic models of $\CT$ and, hence, $I(\CT)\geq \k$.
Now  $\prod _{i\in I}I(\X _i)=\max \{ I(\X _i):i\in I\}=I(\X _{i_1})>\o$, for some $i_1\in I$, by Morley's theorem we have $I(\X _{i_1})\in \{ \o _1,\c\}$
and, as above,  $I(\CT)\geq I(\X _{i_1})$.
\kdok
\begin{te}\label{T709}
If $\CT$ is an FLD$_1$-theory having an atomic model, $\X ^{\rm at}$,
then for each decomposition $\sum _{\BI}\X_i^{\rm at}\in \CD (\X ^{\rm at})$ we have
$$\textstyle
I(\CT)\leq \prod _{i\in I}I(\X _i^{\rm at}).
$$
\end{te}
\dok
Let $\X ^{\rm at}=\sum _{\BI}\X_i^{\rm at}$, where $\BI=\la n,\leq _{\BI}\ra$,
let $r_i^{\rm at}:=\max \X _i^{\rm at}$, for $i<n$, and $\bar r ^{\rm at}:=\la r_0^{\rm at},\dots ,r_{n-1}^{\rm at}\ra$.
We prove first that for each model $\Y$ of $\CT$ we have
\begin{equation}\label{EQ745}\textstyle
\X^{\rm at}\preccurlyeq \Y \;\Rightarrow\; \Y =\sum _{\BI}\Y_i, \mbox{ where }\X^{\rm at} _i\preccurlyeq\Y_i, \mbox{ for each }i<n.
\end{equation}
Let  $\X^{\rm at}\preccurlyeq \Y$.
By Theorem \ref{T700}(a) we have $X_i^{\rm at}:=D_{\f _i (\bar r ^{\rm at} , v), \X ^{\rm at}}$, for $i<n$,
and $X^{\rm at}/D_{\ve (\bar r^{\rm at},u,v),\X^{\rm at} } =\{ X_i^{\rm at}:i<n\}$.
Consequently we have $\X^{\rm at}\models \f _p[\bar r ^{\rm at}]$,
where $\f_p (\bar w)$ is the formula saying that $\{ D_{\f _i (\bar w,v), \cdot}: i<n\}$ is a partition of the domain
and defined by (\ref{EQ714}).
So, since $\X^{\rm at}\preccurlyeq \Y$, we have $\Y\models \f _p[\bar r ^{\rm at}]$
and, defining $Y_i:=D_{\f _i (\bar r ^{\rm at} , v), \Y}$, for $i<n$,
we obtain a partition $\{ Y _i :i<n\}$ of the set $Y$ and
\begin{equation}\label{EQ737}
Y/D_{\ve (\bar r^{\rm at},u,v),\Y} =\{ Y_i:i<n\}= \{ D_{\f _i (\bar r ^{\rm at} , v), \Y}:i<n\}.
\end{equation}
Since $\X^{\rm at} =\sum _{\BI}\X_i^{\rm at}$
we have $\X^{\rm at}\models \f _m[\bar r ^{\rm at}]$,
where $\f_m (\bar w)$ is the formula describing the order between different summands
and defined by (\ref{EQ716}).
Thus $\Y\models \f _m[\bar r ^{\rm at}]$,
since $Y_i:=D_{\f _i (\bar r ^{\rm at} , v), \Y}$, for $i<n$, by (\ref{EQ716}) we have
\begin{equation}\label{EQ738}\textstyle
\forall i,j<n \;\;(i\neq j \Rightarrow \forall y\in Y _i \;\;\forall y'\in Y_j \;\;(y\leq _\Y y' \Leftrightarrow i<_\BI j)),
\end{equation}
and, as in the proof of Theorem \ref{T700}(c), we show that $\Y =\sum _{\BI}\Y_i$.
For $i<n$ we have $r_i^{\rm at}=\max \X _i^{\rm at}$
and, hence, $\X^{\rm at}\models \f _{r_i}[\bar r ^{\rm at}]$,
where $\f_{r_i} (\bar w)$ is the formula defined by (\ref{EQ721}).
Consequently we have $\Y\models \f _{r_i}[\bar r ^{\rm at}]$
and, hence, $r_i^{\rm at}=\max D_{\f _i (\bar r ^{\rm at} , v), \Y}=\max Y_i$.
Finally, for $i<n$ we prove that  $\X^{\rm at} _i\preccurlyeq\Y_i$.
Since $X_i^{\rm at}:=D_{\f _i (\bar r ^{\rm at} , v), \X ^{\rm at}}$, $Y_i:=D_{\f _i (\bar r ^{\rm at} , v), \Y}$ and $\X^{\rm at}\preccurlyeq \Y$
by Fact \ref{T706} we have
$\X^{\rm at} _i=\X^{\rm at} \upharpoonright D_{\f _i (\bar r ^{\rm at} , v), \X ^{\rm at}}\preccurlyeq \Y \upharpoonright D_{\f _i (\bar r ^{\rm at} , v), \Y}=\Y _i$
and (\ref{EQ745}) is proved.

Clearly, for $i<n$ we have
\begin{equation}\label{EQ725}\textstyle
\k _i := I(\X _i^{\rm at})\leq \k:=\prod _{i<n}I(\X _i^{\rm at});
\end{equation}
let $\Mod (\Th (\X^{\rm at} _i),\o)/\!\cong \;=\{ [\Z _i^j]:j<\k _i\}$ be an enumeration.
We show that
\begin{equation}\label{EQ739}\textstyle
\Mod (\CT,\o)/\!\cong \;=\{ [\sum _{\BI}\Z _i^{j_i}]: \la j_0,\dots ,j_{n-1}\ra \in \prod _{i<n}\k _i\}.
\end{equation}
If $[\Y ]\in \Mod (\CT,\o)/\!\cong$,
then w.l.o.g.\ we assume that $\X^{\rm at}\preccurlyeq \Y$
and by (\ref{EQ745}) we have $\Y =\sum _{\BI}\Y_i$, where for $i<n$ we have $\X^{\rm at} _i\preccurlyeq\Y_i$,
and, hence, $\Y_i\cong \Z _i^{j_i}$, for some $j_i<\k _i$.
By Fact \ref{T200}(a) we have $\Y \cong\sum _{\BI}\Z _i^{j_i} $, that is $[\Y ]=[\sum _{\BI}\Z _i^{j_i}]$.
Conversely, if $\la j_0,\dots ,j_{n-1}\ra \in \prod _{i<n}\k _i$,
then for $i<n$ we have $\Z _i^{j_i}\equiv \X^{\rm at} _i$
and, by Fact \ref{T200}(c), $\sum _{\BI}\Z _i^{j_i}\equiv \sum _{\BI} \X^{\rm at} _i=\X^{\rm at}$.
Thus $[\sum _{\BI}\Z _i^{j_i}]\in \Mod (\CT,\o)/\!\cong$
and (\ref{EQ739}) is true.
Now, by (\ref{EQ739}) and (\ref{EQ725}) we have $I(\CT)\leq \prod _{i<n}\k _i=\prod _{i\in I}I(\X _i^{\rm at})$.
\kdok
\begin{te}\label{T717}
If $\CT$ is an FLD$_1$-theory, then we have the following cases:\\[-7mm]
\begin{itemize}\itemsep=-1mm
\item[\sc (i)]   There is $\sum _{\BI}\X_i\in \CD (\CT)$ such that $I(\X _i)=\c$, for some $i\in I$, or $\CT$ is large; then $I(\CT)=\c$;
\item[\sc (ii)]  There is $\sum _{\BI}\X_i\in \CD (\CT)$ such that $I(\X _i)=1$, for all $i\in I$; then $I(\CT)=1$;
\item[\sc (iii)] Otherwise, $\CT$ is small, has an atomic model, $\X ^{\rm at}$,
                 \begin{equation}\label{EQ742}\textstyle
                 \forall \sum _{\BI}\X_i\in \CD (\CT)\;\Big((\forall i\in I \;\;I(\X _i)<\c) \land (\exists i\in I \;\;I(\X _i)>1)\Big),
                 \end{equation}
                 and we have the following subcases:\\[-7mm]
                 \begin{itemize}\itemsep=-1mm
                 \item[\sc (iii.1)] There is $\sum _{\BI}\X_i^{\rm at}\in \CD (\X ^{\rm at})$ such that $\X _i^{\rm at}$ satisfies VC for each $i\in I$, then $3\leq I(\CT)\leq \o$;
                 \item[\sc (iii.2)] For each $\sum _{\BI}\X_i^{\rm at}\in \CD (\X ^{\rm at})$ there is $i\in I$ such that $I(\X _i^{\rm at})=\o_1<\c$; then $I(\CT)=\o _1 <\c$.
                 \end{itemize}
\end{itemize}
\end{te}
\dok
If {\sc (i)} holds and $\CT$ is large, then, clearly, $I(\CT)=\c$.
Otherwise, by Theorem \ref{T711} we have $I(\CT)\geq\c$ and, hence, $I(\CT)=\c$ again.
If {\sc (ii)} holds, then by Theorem \ref{T710} we have $I(\CT)=1$.

Under the assumptions of {\sc (iii.1)}
by (\ref{EQ742})  we have $I(\X^{\rm at} _i)\leq \o$, for all $i\in I$;
so, by Theorem \ref{T709}, $I(\CT)\leq \prod _{i\in I}I(\X^{\rm at} _i)\leq \o$.
By (\ref{EQ742}) again there is $i\in I$ such that $I(\X^{\rm at} _i)>1$
and, by Theorem \ref{T710}, $I(\CT)>1$;
thus, by Vaught's theorem, $I(\CT)\geq 3$.

In subcase {\sc (iii.2)}
for each decomposition $\sum _{\BI}\X ^{\rm at}_i\in \CD (\X ^{\rm at})$ we have:
(a) $I(\X ^{\rm at}_i)= \o _1<\c$, for some $i\in I$;
(b) $I(\X ^{\rm at}_i)\leq \o _1$, for all $i\in I$, (by (\ref{EQ742}) and Morley's theorem).
So, by Theorems \ref{T709} and \ref{T711},  $I(\CT)\leq \prod _{i\in I}I(\X ^{\rm at}_i)=\o _1 \leq I(\CT)$,
which gives $I(\CT)=\o _1 <\c$.
\kdok
If $\CT$ is an FLD$_1$-theory,
then a decomposition $\sum _{\BI}\X_i\in \CD (\CT)$ will be called a {\it VC-decomposition} (resp.\ a {\it VC$^\sharp$-decomposition})
iff $\X _i$ satisfies VC (resp.\ VC$^\sharp$) for each $i\in I$.
\begin{te}\label{T720}
An FLD$_1$-theory $\CT$ satisfies VC iff $\CT$ is large or its atomic model has a VC decomposition.
\end{te}
\dok
If $I(\CT)=\o _1<\c$,
then we have Subcase {\sc (iii.2)} in Theorem \ref{T717};
so $\CT$ is small and $\X ^{\rm at }$ has no VC decomposition.
Conversely, let $\CT$ be a small theory and let $\sum _{\BI}\X_i^{\rm at}$ be a VC decomposition of $\X ^{\rm at }$.
If $\prod _{i<n}I(\X^{\rm at} _i)=\c$, then we have Case {\sc (i)} in Theorem \ref{T717} and $I(\CT)=\c$.
Otherwise we have $\prod _{i<n}I(\X^{\rm at} _i)\leq\o$ and by, Theorem \ref{T709}, $I(\CT)\leq\o$.
\kdok
\begin{rem}\label{R700}\rm
1. The following statements are equivalent (in ZFC):
(a) VC,
(b) VC for complete theories of partial order,
(c) VC for FLD$_1$-theories.
Namely, (a) $\Leftrightarrow$ (b) is well-known (see \cite{Hodg}, p.\ 231) and (b) $\Rightarrow$ (c) is trivial.
If (b) is false, $I(\X)=\o _1<\c$ and $\Y$ is obtained by adding a largest element to $\X$,
then $\Th (\Y)$ is an FLD$_1$-theory and $I(\Y)=\o _1<\c$, thus (c) is false.

2. If the statement ``VC is preserved under finite lexicographic sums  of partial orders with a largest element"
is not a theorem of ZFC and $\X=\sum _{\BI}\X_i$ is a counterexample,
then by Theorem \ref{T717} $\Th (\X)$ is small,
$\prod _{i\in I}I(\X _i)\leq \o$
and $\prod _{j\in J} I(\X ^{\rm at}_j)=\o _1< \c$, for each $\sum _{\BJ}\X ^{\rm at}_j\in \CD (\X ^{\rm at})$.
\end{rem}
\section{Sufficient conditions for VC. Duals }\label{S4}
Regarding condition (iii) in Theorem \ref{T700}(c),
an FLD$_1$-theory $\CT$ will be called {\it actually Vaught's}
iff there are a decomposition $\sum _{\BI}\X_i\in \CD (\CT)$ and sentences $\t _i^{\rm vc} \in \Th (\X _i)$, for $i\in I$,
providing VC; namely,
\begin{equation}\label{EQ741}\textstyle
\exists \sum _{\BI}\X_i\in \CD (\CT) \;\;\forall i\in I \;\;\exists \t _i^{\rm vc} \in \Th (\X _i) \;\;\forall \Z\models \t _i^{\rm vc} \;\;(I(\Z)\leq \o \lor I(\Z)=\c).
\end{equation}
\begin{te}\label{T712}
Vaught's conjecture is true for each actually Vaught's  FLD$_1$-theory $\CT$, more precisely,
\begin{equation}\label{EQ754}\textstyle
I(\CT)   = \left\{
           \begin{array}{cl}
           \c,         & \mbox{ if }\; \exists \sum _{\BI}\X_i\in \CD (\CT)\;\exists i\in I\; I(\X _i)=\c, \;\mbox{ or $\CT$ is large}, \\[1mm]
           1,          & \mbox{ if }\; \exists \sum _{\BI}\X_i\in \CD (\CT)\;\forall i\in I\; I(\X _i)=1,\\[1mm]
           \in [3,\o], & \mbox{ otherwise.}
           \end{array}
           \right.
\end{equation}
\end{te}
\dok
By Theorem \ref{T717} we have to prove that in Case {\sc (iii)} we have Subcase {\sc (iii.1)}.
Namely, assuming that $\X ^{\rm at}$ is an atomic model of $\CT$ and that (\ref{EQ742}) holds,
we prove that there is $\sum _{\BI}\X_i^{\rm at}\in \CD (\X ^{\rm at})$ such that $\X _i^{\rm at}$ satisfies VC, for each $i\in I$.
So, since $\CT$ is actually Vaught's,
there are $\sum _{\BI}\X_i\in \CD (\CT)$, where $I=n\in \N$, and $\t _i^{\rm vc} \in \Th (\X _i)$, for $i<n$,
such that
\begin{equation}\label{EQ756}
\forall i<n\;\; \forall \Z\models \t _i^{\rm vc} \;\;(I(\Z)\leq \o \;\lor\; I (\Z)=\c ).
\end{equation}
By Theorem \ref{T700}(c) there is $\bar r ^{\rm at}:=\la r_0^{\rm at},\dots ,r_{n-1}^{\rm at}\ra \in (X^{\rm at})^n$
such that defining $X_i^{\rm at}:=D_{\f _i (\bar r ^{\rm at} , v), \X ^{\rm at}}$, for $i<n$, we have:
(i) $\{ X_i^{\rm at}:i<n\}$ is a partition of the set $X^{\rm at}$ and $X^{\rm at}/D_{\ve (\bar r^{\rm at},u,v),\X^{\rm at} } =\{ X_i^{\rm at}:i<n\}$;
(ii) $\X^{\rm at} =\sum _{\BI}\X_i^{\rm at}$ and $r_i^{\rm at}=\max \X _i^{\rm at}$, for $i<n$;
(iii) $\X _i^{\rm at} \models \t _i^{\rm vc}$, for $i<n$.
Now, by (ii) we have $\sum _{\BI}\X_i^{\rm at}\in \CD (\X ^{\rm at})$
and by (iii), (\ref{EQ756}) and (\ref{EQ742}) we have $I(\X _i^{\rm at})\leq \o$, for all $i<n$.
\kdok
\begin{te}\label{T713}
If $\CT$ is an FLD$_1$-theory having a {\rm VC}$^{\sharp}$ decomposition, then {\rm VC}$^{\sharp}$ holds for $\CT$.
\end{te}
\dok
If $\sum_{\BI}\X_i\in \CD (\CT)$ is a {\rm VC}$^{\sharp}$ decomposition, then $I(\X_i)\in \{1,\c\}$, for all $i\in I$.
So, in Theorem \ref{T717} we have Case {\sc (i)} or Case {\sc (ii)}, which gives $I(\CT)\in \{1,\c\}$.
\kdok
\paragraph{Duals of Theorems \ref{T700}--\ref{T713}}
Dually we define partial orders admitting a {\it finite lexicographic decomposition with zeros},
FLD$_0$-{\it partial orders} (in a decomposition $\X =\sum _{\BI}\X_i$ we require that $\min \X _i$ exists, for each $i\in I$).
Then $\I (\X):=\{ \BI \in \CC ^{\rm fin}: \mbox{ there is an FLD$_0$ decomposition } \X =\sum _{\BI}\X_i\}$,
a complete theory of partial order $\CT$ is called an FLD$_0$-theory
iff some (equivalently, each) model of $\CT$ is an FLD$_0$-partial order,
and we define $\I(\CT):= \I (\X)$, where $\X$ is some (any) model of $\CT$.
Finally, $\CD (\CT)$ is the class of all FLD$_0$-decompositions of models of $\CT$ defined by (\ref{EQ740})
and $\CT$ is an {\it actually Vaught's} FLD$_0$-theory iff (\ref{EQ741}) holds.
Thus,  writing FLD$_0$ instead of FLD$_1$ in Theorems \ref{T700}--\ref{T713} we obtain their duals.
We will not list them explicitly.
\section{Closures satisfying VC}\label{S5}
Let $\CC$ be an $\cong$-closed class of partial orders
and let $\la \CC \ra _{\Sigma}$ be the minimal class of partial orders containing $\CC$
and which is closed under isomorphism and finite lexicographic sums.\footnote{If $\BI \in \CC ^{\rm fin}$
and $\la \Y _i :i\in I\ra \in \CC ^I$,
then in $\CC$ there are $\X_i \cong\Y _i $, for $i\in I$, with pairwise disjoint domains
and  $\sum _\BI \X _i\in \la \CC \ra _{\Sigma}$.
Taking another representatives $\X_i' \cong\Y _i $
by Fact \ref{T200}(a) we have $\sum _\BI \X _i'\cong \sum _\BI \X _i$;
so $\la \CC \ra _{\Sigma}$ can be regarded as a closure under finite lexicographic sums of order types from $\CC$.}
For a description of $\la \CC \ra _{\Sigma}$ we first prove that a lexicographic sum of lexicographic sums is a lexicographic sum.
\begin{fac}\label{T541}
If $\BI$, $\BJ_i$, $i\in I$, and $\X _j $, $j\in \bigcup _{i\in I}J_i$, are partial orders with pairwise disjoint domains, then
\begin{equation} \label{EQ621}\textstyle
\sum _{\BI}\sum _{\BJ _i}\X _j=\sum _{\sum _{\BI}\BJ _i }\X _j.
\end{equation}
\end{fac}
\dok
Let $\BI=\la I, \leq _{\BI}\ra$, $\BJ_i=\la J_i, \leq _{\BJ_i}\ra$, for $i\in I$, $\X _j =\la X_j, \leq _{\X _j}\ra$, for $j\in \bigcup _{i\in I}J_i$,
and $\BJ =\sum _{\BI}\BJ _i=\la J,\leq _{\BJ}\ra$.
Clearly $J=\bigcup _{i\in I}J_i$ and the posets from (\ref{EQ621}) have the same domain: $X=\bigcup _{j\in J}X_i$.
Let $\sum _{\BI}\sum _{\BJ _i}\X _j=\la X, \leq \ra$,  $\sum _{\BJ }\X _j=\la X, \leq '\ra$ and $x,y\in X$.

Assuming that $x\leq y$ we prove that $x\leq' y$.
First, let $x,y\in \bigcup _{j\in J_i}\X _j$, for some $i\in I$, and $x\leq _{\sum _{\BJ _i}\X _j}y$.
If for some $j\in J_i$ we have $x,y\in X_j$ and $x\leq _{\X _j}y$, then, clearly, $x\leq 'y$.
Otherwise, there are different $j,j'\in J_i$ such that $x\in X_j$, $y\in X_{j'}$ and $j < _{\BJ _i}j'$; then $j < _{\BJ }j'$ and, hence, $x\leq 'y$ again.
Second, if $x\in \bigcup_{j\in J_i} X_j$ and $y\in \bigcup_{j\in J_{i'}} X_j$, where $i<_{\BI}i'$,
then $x\in X_j$, for some $j\in J_i$, $y\in X_{j'}$, for some $j'\in J_{i'}$ and, since $i<_{\BI}i'$, we have $j<_{\BJ}j'$;
thus, $x\leq 'y$ indeed.

Conversely, assuming that $x\leq ' y$ we prove that $x\leq y$.
First, if there are $i\in I$ and $j\in J_i$ such that $x,y\in X_j$ and $x \leq_{\X _j} y$, then $x\leq _{\sum _{\BJ _i}\X _j}y$ and, hence, $x\leq y$.
Second, let $x\in X_j$ and $y\in X_{j'}$, where $j<_{\BJ} j'$.
If $j,j'\in J_i$, for some $i\in I$, and $j<_{\BJ_i} j'$, then $x\leq _{\sum _{\BJ _i}\X _j}y$ and, hence, $x\leq y$.
Otherwise, there are different $i,i'\in I$ such that $j\in J_i$ and $j'\in J_{i'}$.
Then, since $j<_{\BJ} j'$ we have $i<_{\BI}i'$
and since $x\in \bigcup_{j\in J_i}X_j$ and $y\in \bigcup_{j\in J_{i'}}X_j$
we have $x\leq y$.
\kdok
\begin{fac}\label{T722}
If $\CC$ is a class of partial orders closed under isomorphism, then
\begin{equation}\label{EQ757}\textstyle
\la \CC \ra _{\Sigma}=\bigcup\{\IT (\sum _\BI \X _i):\BI \in \CC ^{\rm fin}\land  \la \X _i :i\in I\ra \in \CC ^I \land \forall i,j\in I \; (i\neq j\Rightarrow X_i \cap X_j=\emptyset)\}.
\end{equation}
\end{fac}
\dok
Let $\CC^*$ denote the r.h.s.\ of (\ref{EQ757}).
First, for $\Y\in \CC$ we have $\Y =\sum _1 \Y\in \CC ^*$; thus $\CC \subset \CC ^*$.
Second, it is evident that $\CC ^*$ is $\cong$-closed.
Third, the class $\CC ^*$ is closed under finite lexicographic sums
because by Fact \ref{T541} a lexicographic sum of lexicographic sums of elements of $\CC $ is a lexicographic sum of elements of $\CC $.
Finally, if a class $\CC' \supset \CC$ is closed under $\cong$ and finite lexicographic sums, then, clearly, $\CC ^* \subset \CC '$.
\kdok
We recall that a partial order $\X$ is a (model-theoretic) {\it tree} iff $(\cdot ,x]$ is a linear order, for each $x\in X$,
and that $\X$ is a {\it reticle} iff it does not embed the four-element poset with the Hasse diagram $N$.
Note that adding a smallest (or a largest) element to a reticle produces a reticle again.
In \cite{Sch3} Schmerl confirmed VC for reticles and proved that the theory of reticles is finitely axiomatizable; see Corollary 4.7 of \cite{Sch3}.
Thus the classes
$$
\CC ^{\rm lo}_0\subset\CC ^{\rm tree}_0\subset \CC ^{\rm reticle}_0 \mbox{ and }\CC ^{\rm ba}
$$
of linear orders with a smallest element, rooted trees, reticles with a smallest element and Boolean algebras
are first-order definable by the sentences $\bigwedge\CT ^{\rm lo}_0$, $\bigwedge\CT ^{\rm tree}_0$, $\bigwedge\CT ^{\rm reticle}_0$ and $\bigwedge\CT ^{\rm ba}$.

Let $\CC ^{\rm fin}_0$ (resp.\ $\CC ^{\rm fin}_1$) denote the class of finite partial orders with a smallest (resp.\ largest) element.
If $\CC$ is a class of partial orders,
by $\CC ^{-1}$ we denote the class of the corresponding reversed orders, $\X ^{-1}:=\la X, (\leq _\X)^{-1}\ra$, for $\X\in \CC$.
So, $(\CC ^{\rm tree}_0)^{-1}$ is the class of reversed trees with a largest element,
$(\CC ^{\rm fin}_0)^{-1}=\CC ^{\rm fin}_1$,
$(\CC ^{\rm reticle}_0)^{-1}=\CC ^{\rm reticle}_1$ (the class of reticles with a largest element) and
$(\CC ^{\rm ba})^{-1}=\CC ^{\rm ba}$.
\begin{te}\label{T714}
Vaught's conjecture (in fact (\ref{EQ754})) is true for the theory of each partial order $\sum _{\BI}\X_i$ from the class
$$
 \Big\la \CC ^{\rm fin}_0 \cup \CC ^{\rm reticle}_0\cup \CC ^{\rm ba}\Big\ra_{\Sigma}
\cup \Big\la \CC ^{\rm fin}_1 \cup \CC ^{\rm reticle}_1\cup \CC ^{\rm ba}\Big\ra_{\Sigma}.
$$
In particular, Vaught's conjecture is true for lexicographic sums of rooted trees, Boolean algebras etc.
\end{te}
\dok
Let $\X:=\sum _{\BI}\X_i \in \la \CC ^{\rm fin}_0 \cup \CC ^{\rm reticle}_0\cup \CC ^{\rm ba}\ra_{\Sigma}$.
In order to apply the dual of Theorem \ref{T712}
we note that $\Th (\X)$ is an FLD$_0$-theory and show that it is actually Vaught's.
So, if $i\in I$ and $\X _i \in \CC ^{\rm reticle}_0$, that is, if $\X_i$ is a reticle with a smallest element,
then $\t _i:=\bigwedge\CT ^{\rm reticle}_0 \in \Th (\X _i)$
and, by Schmerl's result, VC is true for each $\Z\models \t _i$.
If $\X _i \in \CC ^{\rm ba}$ the same holds by the result of Iverson \cite{Ive}
and if $\X _i \in \CC ^{\rm fin}_0$, we have a triviality.
Thus $\Th (\X)$ is actually Vaught's.
For $\X\in \la \CC ^{\rm fin}_1 \cup \CC ^{\rm reticle}_1\cup \CC ^{\rm ba}\ra_{\Sigma}$
we have a dual proof.
\kdok
In order to extend the result of Fact \ref{T719} (concerning disconnected partial orders) we introduce a new closure.
First, if $\X =\la X,\leq _\X \ra$ is a partial order,
let $1+\X $ be the partial order obtained from $\X$ by adding an element, say $x_0 \not\in X$, below all elements of $X$
(thus, $1+\X =\la X \cup \{ x_0\},\leq _{1+\X}\ra$,
where  $\leq _{1+\X}\;=\;\leq_\X \;\cup\;\{ \la x_0 ,x\ra :x\in X \cup \{ x_0\}\}$).
Second, for any partial order $\X$ let us define the {\it rooted $\X$}, $\X _r$, by
\begin{equation}\label{EQ749}\textstyle
\X _r    = \left\{
           \begin{array}{cl}
           \X,        & \mbox{ if } \min \X \mbox{ exists},  \\
           1+\X,         & \mbox{ otherwise}.
           \end{array}
           \right.
\end{equation}
Third, if $\BI \in \CC^{\rm fin}$ and $\X _i$, $i\in I$, are partial orders,
let $\sum _{\BI}^r \X _i := \sum _{\BI} (\X _i)_r$ be the corresponding {\it lexicographic sum of rooted summands} $\X _i$.
Now, for an $\cong$-closed class $\CC$ of partial orders
let $\la \CC\ra _{\Sigma ^r}$ be the minimal closure of $\CC$ under isomorphism and finite lexicographic sums of rooted summands.
Clearly, $\la \CC\ra _{\Sigma ^r}= \bigcup _{n\in \o} \CC _n$, where $\CC _0:=\CC$ and, for $n\in \o$,
\begin{equation}\label{EQ750}\textstyle
\CC _{n+1}:=\bigcup \{ \IT(\sum _{\BI} (\X _i)_r) : \BI \in \CC^{\rm fin} \land \la \X _i:i\in I\ra \in (\bigcup _{m\leq n}\CC _m)^I
                                          \land \forall \{i,j\}\in [I]^2 \; (X_i)_r \cap (X_j)_r=\emptyset\}.
\end{equation}
For example, for $n=1$, $\sum _{\BI} (\X _i)_r=\sum _{\BI} (\sum _{\BJ _i} (\X _i^j)_r)_r \in \CC_2$, where $\X _i^j \in \CC$, for $i\in I$ and $j\in J_i$.
Let $\CC ^{{\rm VC}^{\sharp}}$ be the class of \emph{all} partial orders satisfying {\rm VC}$^{\sharp}$.
\begin{te}\label{T715}
If $\CC$ is an $\cong$-closed class of partial orders satisfying {\rm VC}$^{\sharp}$,
then {\rm VC}$^{\sharp}$ holds for each partial order from the closure $\la \CC\ra _{\Sigma ^r}$.
In particular, the class $\CC ^{{\rm VC}^{\sharp}}$
is closed under under finite lexicographic sums of rooted summands,
that is $\la \CC ^{{\rm VC}^{\sharp}}\ra _{\Sigma ^r}=\CC ^{{\rm VC}^{\sharp}}$.
\end{te}
\dok
By induction we prove that for each $n\in \o$ each partial order $\X\in \CC_n$ satisfies {\rm VC}$^{\sharp}$.
For $n=0$ this is our hypothesis.
Let the statement be true for all $m\leq n$ and let $\X=\sum _{\BI} (\X _i)_r\in \CC _{n+1}$.
Then for each $i\in I$ the partial order $(\X _i)_r$ has a smallest element,
so, $\X$ is an FLD$_0$-poset and, by the dual of Theorem \ref{T708}, $\Th (\X)$ is an FLD$_0$-theory.
By (\ref{EQ750}) for each $i\in I$ we have $\X _i \in \bigcup _{m\leq n}\CC _m$,
by the induction hypothesis the poset $\X _i$ satisfies {\rm VC}$^{\sharp}$
and, hence, the poset $(\X _i)_r$ satisfies {\rm VC}$^{\sharp}$ too.
So, by the dual of Theorem \ref{T713}, $\X $ satisfies {\rm VC}$^{\sharp}$ as well.
\kdok
Theorem \ref{T715} provides the following extension of Fact \ref{T719}.
Recall that $\CC ^{\rm lo}$, $\CC ^{\rm tree}_{\rm 0,fmd}$ and $\CC ^{\rm tree}_{\rm if,{\rm VC}^{\sharp}}$
are the classes of linear orders, rooted FMD trees, and initially finite trees satisfying  {\rm VC}$^{\sharp}$
and that $\la \CC ^{\rm lo}\ra_{\dot{\cup}_{\infty}}$ is the class of infinite disjoint unions of linear orders.
\begin{te}\label{T716}
{\rm VC}$^{\sharp}$ is true for the theory of each partial order from the class
\begin{equation}\label{EQ751}\textstyle
\Big\la \CC ^{\rm fin} \cup
\la \CC ^{\rm lo} \ra _{\dot{\cup}\Pi} \cup
\la \CC ^{\rm ba}\ra _{\dot{\cup}\Pi} \cup
\la \CC ^{\rm tree}_{\rm 0,fmd}\ra _{\dot{\cup}\Pi} \cup
\la \CC ^{\rm tree}_{\rm if,{\rm VC}^{\sharp}}\ra _{\dot{\cup}\Pi}\cup
\la \CC ^{\rm lo}\ra_{\dot{\cup}_{\infty}}\Big\ra_{\Sigma ^r}.
\end{equation}
In particular, {\rm VC}$^{\sharp}$ is true for finite lexicographic sums of finite products of linear orders with zero, Boolean algebras, rooted FMD trees etc.
\end{te}
Theorem \ref{T716} generates a jungle of partial orders satisfying {\rm VC}$^{\sharp}$.
Namely, if $\CC$ is an $\cong$-closed class of partial orders,
then by Lemma 3.1 of \cite{Ksharp}  $\X \in \la \CC \ra _{\dot{\cup}\Pi}$ iff $\X=\dot{\bigcup}_{i<n}\prod _{j<m_i}\X _i^j$, for some $n,m_i\in \N$ and $\X _i^j\in \CC$;
so, {\rm VC}$^{\sharp}$ is true for partial orders of the form
$$\textstyle
\sum _{\BI}(\dot{\bigcup}_{j<n_i}\prod _{k<m_i^j}\X _i^{j,k})_r,
$$
where for each $i\in I$ we have: $\X _i^{j,k}\in \CC$,
where $\CC \in \{ \CC ^{\rm lo},\CC ^{\rm tree}_{\rm 0,fmd},\CC ^{\rm tree}_{\rm if,{\rm VC}^{\sharp}},\CC ^{\rm ba}\}$, for all $j<n_i$ and $k<m_i^j$.

Theorem \ref{T716} and the operation $\sum ^r$ are related to FLD$_0$-posets
and in a natural way we obtain a dual statement and operation related to FLD$_1$-posets; e.g., in (\ref{EQ749}), instead of $1+\X$ we take $\X +1$ etc.
\begin{rem}\rm
{\it The closures $\la\CC \ra _{\dot{\cup}\Pi}$, $\la\CC \ra _{\Sigma}$ and $\la\CC \ra _{\Sigma ^r}$.}
If $\CC$ is a $\cong$-closed class of posets,
then by Lemma 3.1 of \cite{Ksharp} and Fact \ref{T722}
its closures $\la\CC \ra _{\dot{\cup}\Pi}$ and $\la\CC \ra _{\Sigma}$ are obtained in one step.
Concerning the closure $\la\CC \ra _{\Sigma ^r}$ the situation is different and depends of $\CC$.
For example, by Theorem \ref{T715} for $\CC= \CC^{{\rm VC}^{\sharp}}$ we have $\la \CC \ra _{\Sigma ^r}=\CC $;
so, we do not obtain new structures in the closure.
If we take $\CC =\{ 1\}$, more precisely, if $\CC$ is the class of all one-element posets,
then by (\ref{EQ750}) $\CC _1:=\{ \sum _{\BI} 1 : \BI \in \CC^{\rm fin}\}=\CC ^{\rm fin}$
and since $|\sum _{\BI} (\X _i)_r|<\o$, if $|X_i|<\o$, for all $i\in I$,
by (\ref{EQ750}) we have $\la\CC ^{\rm fin}\ra _{\Sigma ^r}=\CC ^{\rm fin}$,
which implies that $\la\CC \ra _{\Sigma ^r}=\CC _1$,
that is, the closure of $\CC$ is obtained in the first step of the recursion.

Generally, the class $\CC _{n+1}$ defined by (\ref{EQ750}) can be obtained from $\CC _n$ in two steps:
\begin{eqnarray*}
\CC _{n+1}' &  = &\textstyle \bigcup \{ \IT(1+\X) : \X \in \bigcup _{m\leq n}\CC _m \land \min \X \mbox{ does not exist}\} \mbox{ and } \\
\CC _{n+1}  &  = &\textstyle \CC _{n+1}' \cup \la \{ \X \in \CC _{n+1}': \min \X \mbox{ exists}\}\ra _{\Sigma}.
\end{eqnarray*}
So, for the class $\CC$ defined by (\ref{EQ751}) we will have $1+\X \in \CC _1'$, whenever $\X$ is a disjoint union of more than one poset
or, for example, if $\X$ is a direct product of linear orders without a smallest element (e.g.\ $\Z \times \o$)
and in $\CC _2$ we will have all finite lexicographic sums of these ``rooted" posets.
But this is not the end;
namely, for some $\BI$ the posets $\sum _{\BI} (\X _i)_r $ from $\CC _{n+1}$ are without a smallest element;
for example if $\BI$ is an antichain of size $>1$ and we obtain new (isomorphism types of) posets.
\end{rem}

{\footnotesize

\end{document}